\documentclass[12pt,reqno]{amsart}
\ifx\pdfsuppressptexinfo\undefined\else\pdfsuppressptexinfo=-1\fi
\usepackage{amsmath,amsthm,amssymb}
\usepackage[hmargin=1.30in,vmargin=1.04in]{geometry}
\usepackage[colorlinks,linkcolor=blue,citecolor=blue,urlcolor=blue]{hyperref}
\usepackage{bookmark}
\numberwithin{equation}{section}
\newcommand{\F}{\mathbb F}
\DeclareMathOperator{\PGU}{PGU}
\newtheorem{theorem}{Theorem}[section]
\newtheorem{corollary}[theorem]{Corollary}
\newtheorem{proposition}[theorem]{Proposition}
\newtheorem{lemma}[theorem]{Lemma}
\theoremstyle{definition}

\theoremstyle{remark}
\newtheorem{remark}[theorem]{Remark}
\hypersetup{pdftitle={Explicit plane models for all Galois quotients of the Hermitian curve},pdfauthor={Liming Ma and Yipeng Wang}}
\begin{document}
\title[Galois quotients of the Hermitian curve]
{Explicit plane models for all Galois quotients of the Hermitian curve}
\author[Liming Ma]{Liming Ma}
\address{School of Mathematical Sciences, University of Science and Technology of China,
Hefei, Anhui 230026, China}
\email{lmma20@ustc.edu.cn}
\author[Yipeng Wang]{Yipeng Wang}
\address{School of Mathematical Sciences, University of Science and Technology of China,
Hefei, Anhui 230026, China}
\email{jisunxingfu@mail.ustc.edu.cn}
\date{November 28, 2026}

\begin{abstract}
Let $H=\mathbb F_{q^2}(x,y)$ be the Hermitian function field defined by
$y^q+y=x^{q+1}$. For every subgroup $G$ of its $\mathbb F_{q^2}$-automorphism
group $\PGU(3,q)$, we construct an explicit quotient map and compute a plane
equation for $H^G$ from matrix generators of $G$.
Let $S$ be the stabilizer in $G$ of the common pole $P_\infty$ of $x$ and $y$.
For each coordinate, we multiply the distinct elements of its $S$-orbit
and take the trace of the product from $H^S$ to $H^G$. The resulting
functions $U,V$ generate $H^G$. Their only pole is the place below $P_\infty$,
where they have coprime orders $a=q/|G_x|$ and $b=(q+1)/|G_y|$.
Here $G_x$ and $G_y$ are the stabilizers of $x$ and $y$ in $G$, respectively.
The plane equation relating $U$ and $V$ has total degree $\max(a,b)\leq q+1$.
We recover its coefficients from Laurent expansions at $P_\infty$ by
solving a triangular system of $(a+1)(b+1)/2$ linear equations.
The total delta invariant of the resulting affine model is
$g(H^{G_xG_y})-g(H^G)$, where $g(H^{G_xG_y})=(a-1)(b-1)/2$.
For finite groups acting on function fields over perfect fields, the same
construction gives fixed-field generators from two functions whose only pole
is a common rational place and whose pole orders are coprime.
\end{abstract}
\maketitle
\raggedbottom

\section{Introduction}
Let $q$ be a prime power and put $k=\F_{q^2}$.
An algebraic function field $E/k$ of genus $g(E)$ with full constant
field $k$ is called \emph{maximal} if it has $q^2+1+2qg(E)$ rational places,
attaining the Hasse--Weil bound. The Hermitian function field
$H=k(x,y)$, where $y^q+y=x^{q+1}$, is maximal of genus
$q(q-1)/2$. Its group of $k$-automorphisms is $\PGU(3,q)$, of order
$q^3(q^3+1)(q^2-1)$ \cite[Section~2]{GarciaStichtenothXing2000}.
By a result of Serre, every fixed field $H^G$ with
$G\leq\PGU(3,q)$ is also maximal \cite{Lachaud1987}; see also
\cite[Theorem~10.2]{HirschfeldKorchmarosTorres2008}.
Writing $\mathcal H_q$ for the smooth projective curve with function
field $H$, we call the curves $\mathcal H_q/G$ with function fields
$H^G$ the \emph{Galois quotients} of $\mathcal H_q$.

Two natural questions concern the genus of $H^G$ and an explicit
equation for it. Garcia, Stichtenoth and Xing computed many such
genera via the Hurwitz genus formula and gave generators and
equations in several cases \cite{GarciaStichtenothXing2000}.
Further genus calculations were obtained by Abd\'on and Quoos
\cite{AbdonQuoos2004}.
In odd characteristic, Bassa, Ma, Xing and Yeo determined the
genera arising from subgroups of the stabilizer in $\PGU(3,q)$ of
a rational place \cite{BassaMaXingYeo2013}.
Ma and Xing computed genera of fixed fields for subgroups involving
the involution $(x,y)\mapsto(x/y,1/y)$ \cite{MaXing2019}, while
Montanucci and Zini treated groups preserving neither a point nor
a triangle in the projective plane \cite{MontanucciZiniNoFixed2020}.
For $q\equiv1\pmod4$, Montanucci and Zini determined the complete
spectrum of genera of Galois quotients of $\mathcal H_q$
\cite{MontanucciZini2020}.

Explicit plane equations are also known for several families of quotients.
For odd $q$, coordinatewise squaring on the Fermat model of $\mathcal H_q$
gives the Fermat curve of degree $(q+1)/2$ as a Klein four quotient
\cite[Section~4]{CossidenteHirschfeldKorchmarosTorres2000}.
Cossidente, Korchm\'aros and Torres gave plane models for quotients by
subgroups of prime order \cite{CossidenteKorchmarosTorres2000}. Equations for
other families appear in
\cite{DionigiGatti2025,GarciaKawakitaMiura2006,GattiKorchmaros2024},
and \cite{GattiGhiandoniKorchmaros2025} gives an explicit generator
of the rational fixed field of the full group $\PGU(3,q)$.

Using fixed Hermitian coordinates $x,y$, we construct a plane model and
quotient map for every $G\leq\PGU(3,q)$, including groups fixing no rational
place. We express the total delta invariant of the affine model as a
difference of genera and determine when this model is smooth.

Fix such a subgroup $G$ and put $F=H^G$. The pole divisors of $x,y$ in $H$
are $qP_\infty$ and $(q+1)P_\infty$, respectively, for a rational place
$P_\infty$. Write $Q_\infty=P_\infty\cap F$.
Let $S=\{g\in G:g(P_\infty)=P_\infty\}$ be the decomposition group
of $P_\infty$ in $H/F$, and let
$G_t=\{g\in G:g(t)=t\}$ for $t=x,y$.
Writing $S\cdot t=\{s(t):s\in S\}$, define
\begin{equation}\label{hhh:orbit-generators}
 A_t=\prod_{f\in S\cdot t} f\quad(t=x,y),\qquad
 U=\sum_{gS\in G/S}g(A_x),\qquad V=\sum_{gS\in G/S}g(A_y).
\end{equation}
Since $P_\infty$ is the only pole of $t$, we have $G_t\leq S$.
The orbit product is the norm $A_t=\operatorname N_{H^{G_t}/H^S}(t)$,
and $U,V$ are the traces of $A_x,A_y$ from $H^S$ to $F$.
In $H$, each $g(A_t)$ has its only pole at $gP_\infty$. Distinct
cosets give distinct places, so no pole cancels in either trace.

\begin{theorem}\label{hhh:main}
The numbers $a=q/|G_x|$ and $b=(q+1)/|G_y|$ are coprime positive integers.
The functions $U,V$ in \eqref{hhh:orbit-generators} satisfy the following.
\begin{enumerate}
\item We have $F=k(U,V)$, and the pole divisors in $F$ are
\[
 (U)_\infty=aQ_\infty,\qquad (V)_\infty=bQ_\infty.
\]
\item There is a unique absolutely irreducible polynomial $\Phi\in k[T_1,T_2]$,
monic in $T_2$, such that $\Phi(U,V)=0$. Its total degree is
$\max(a,b)\leq q+1$. With respect to the local parameter $\pi=x/y$ at
$P_\infty$, the expansions of $U,V$ modulo $\pi^{|S|ab}k[[\pi]]$
determine the coefficients of $\Phi$ by a triangular system of
$(a+1)(b+1)/2$ linear equations.
\end{enumerate}
\end{theorem}

The map $(U,V)\colon\mathcal H_q\dashrightarrow\{\Phi=0\}$ has degree
$|G|$ and induces a birational map from $\mathcal H_q/G$ to the plane model.

Since $F$ is maximal, $q$ and $q+1$ are non-gaps at every rational
place of $F$ \cite[Proposition~1.5(iv)]{FuhrmannGarciaTorres1997}.
Thus Lemma~\ref{hhh:coprime-plane} already gives a plane model of degree $q+1$.
Formula~\eqref{hhh:orbit-generators} gives explicit generators of pole orders
$a\mid q$ and $b\mid q+1$. The resulting degree is less than $q+1$ exactly when
$|G_y|>1$. The sections $1,U,V\in\mathcal L((q+1)Q_\infty)$ give the plane map,
and the complete series $|(q+1)Q_\infty|$ gives the natural embedding of
$\mathcal H_q/G$ \cite[Theorem~0.1]{KorchmarosTorres2001}.

Let $R=k[U,V]$ and let $\widetilde R$ be its integral closure in $F$.
The total delta invariant of the affine model is
$\delta_{\mathrm{aff}}=\dim_k(\widetilde R/R)$.
By \eqref{hhh:point-stabilizer}, $G_x$ acts by translations of $y$ and
$G_y$ by scalings of $x$, so $G_xG_y$ is an abelian subgroup.

\begin{samepage}
\begin{theorem}\label{hhh:defect}
Let $D=G_xG_y$. The affine model defined by $U,V$ satisfies
\begin{equation}\label{hhh:defect-eq}
 \delta_{\mathrm{aff}}
 =\frac{(a-1)(b-1)}2-g(F)
 =g(H^D)-g(F).
\end{equation}
Thus $g(F)\le\tfrac12(a-1)(b-1)$, with equality if and only if
the affine model is smooth. If $g(F)\ge2$, these conditions are
equivalent to $G=D$.
\end{theorem}
\end{samepage}

The Weierstrass semigroup at $P_\infty\cap H^D$ is $\langle a,b\rangle$,
and $\delta_{\mathrm{aff}}$ counts the non-gaps at $Q_\infty$ outside this semigroup.
Since $D$ is abelian, Theorem~\ref{hhh:defect} gives a singular model for every
nonabelian $G$ with $g(F)\ge2$.
For $G$ fixing $P_\infty$, \cite[Corollary~5.7]{MaWangDecomposition}
gives generators for a smooth affine plane model. In this case the model
defined by \eqref{hhh:orbit-generators} has coordinates $A_x,A_y$ and is
singular whenever $g(F)\ge2$ and $G\ne D$.

Sections~\ref{sec:general} and~\ref{hhh:geometric} prove
Theorems~\ref{hhh:main} and~\ref{hhh:defect}, respectively.
For finite group actions over perfect fields, Remark~\ref{hhh:general-fields}
constructs fixed-field generators from two functions whose only pole is
the same rational place and whose pole orders are coprime. When these
orders $1<m<n$ generate the Weierstrass semigroup there,
Proposition~\ref{hhh:semigroup} extends \eqref{hhh:defect-eq} and the smoothness criterion.
For $q=2^h$, Section~\ref{sec:s3} treats the group $S_3$ permuting the
Fermat coordinates. It fixes no rational place, and each trace in
\eqref{hhh:orbit-generators} has three summands in the chosen Hermitian
coordinates. Proposition~\ref{hhh:s3-result} gives the quotient equation
and, for $h\ge2$, recovers the genus from the local delta invariant of
its unique affine singularity.

\section{Generators from a decomposition group and its cosets}
\label{sec:general}

We identify places of $H$ with closed points of $\mathcal H_q$. Write
$v_P$ and $\mathcal O_P$ for the normalized valuation and valuation ring
at $P$ \cite{Stichtenoth2009}, and $(f)^E_\infty$ for the pole divisor
of $f$ in $E$, omitting the superscript when the field is clear.

Let $\mathcal H_q$ be given in projective coordinates $(X:Y:Z)$ by
$X^{q+1}=Y^qZ+YZ^q$, so that $x=X/Z$, $y=Y/Z$ and
$P_\infty=(0:1:0)$. By \cite[Proposition~6.4.1(g),(h)]{Stichtenoth2009},
the pole divisors of $x,y$ are $qP_\infty,(q+1)P_\infty$, respectively,
and the Weierstrass semigroup at $P_\infty$ is
$\Gamma(P_\infty)=\langle q,q+1\rangle$.
Thus $\pi=X/Y=x/y$ has valuation $1$ and is a local parameter.
With $\eta=Z/Y=1/y$, the local equation is $\eta+\eta^q=\pi^{q+1}$.

For $g\in\PGU(3,q)$ and $f\in H$, we put $g(f)=f\circ g^{-1}$,
so that $v_{gP}(g(f))=v_P(f)$. The stabilizer of $P_\infty$
in $\PGU(3,q)$ consists of the substitutions
\begin{equation}\label{hhh:point-stabilizer}
 x\longmapsto\alpha x+\beta,\qquad
 y\longmapsto\alpha^{q+1}y+\alpha\beta^q x+\gamma,
\end{equation}
where $\alpha\in k^*$ and $\beta,\gamma\in k$ satisfy
$\gamma^q+\gamma=\beta^{q+1}$
\cite[Section~2]{GarciaStichtenothXing2000}.
Since $P_\infty$ is the only pole of both $x$ and $y$, we have
$G_x,G_y\leq S$.
An element of $G_x$ therefore translates $y$ by some $c\in k$ with $c^q+c=0$.
An element of $G_y$ sends $x$ to $\varepsilon x$ with
$\varepsilon\in k^*$ and $\varepsilon^{q+1}=1$.
In particular $|G_x|\mid q$ and $|G_y|\mid q+1$.

We call $x',y'\in H$ a pair of \emph{Hermitian coordinates} if
$H=k(x',y')$ and $y'^q+y'=x'^{q+1}$.
Such a pair defines a $k$-automorphism $\lambda$ by $(x,y)\mapsto(x',y')$.
Each new coordinate has its only pole at $\lambda(P_\infty)$. Every rational place arises
this way by transitivity of $\PGU(3,q)$. Applying $\lambda$ to the results
for $\lambda^{-1}G\lambda$ in $x,y$ gives them for $G$ in $x',y'$.

\subsection{Pole divisors of the generators}
In our general statements the full constant field $k$ of $E$ may be any
perfect field. We use the following lemma for $x,y$ and again in Proposition~\ref{hhh:semigroup}.

\begin{lemma}\label{hhh:orbit-trace}
Let $E/k$ be a function field with perfect full constant field $k$,
and let $G\le\operatorname{Aut}_k(E)$ be finite. For a rational place
$P$ of $E$, put $S=\operatorname{Stab}_G(P)$ and $Q=P\cap E^G$.
Suppose $t\in E$ has pole divisor $mP$, with $m>0$, and put
$G_t=\{g\in G:g(t)=t\}$ and $A_t=\prod_{f\in S\cdot t} f$.
Then $G_t\le S$, $|G_t|\mid m$, and
\[
 \bigl(\operatorname{Tr}_{E^S/E^G}(A_t)\bigr)_\infty^{E^G}
 =\frac{m}{|G_t|}\,Q.
\]
\end{lemma}
\begin{proof}
Every element of $G_t$ fixes the only pole $P$ of $t$.
The residue degree of $P$ over $Q$ is one, so the ramification index
is $e=|S|$ by \cite[Theorem~3.8.2(a)]{Stichtenoth2009}.
For $h\in E^S$ with $(h)_\infty^E=\ell P$, where $\ell>0$, we first
show that its trace $z$ has pole divisor $(\ell/e)Q$ in $E^G$.
Indeed, $z=\operatorname{Tr}_{E^S/E^G}(h)=\sum_{gS\in G/S}g(h)$.
The cosets parametrize both the $E^G$-embeddings of $E^S$ into $E$
and the places $gP$ above $Q$. Exactly one summand has a pole at each
$gP$, of order $\ell$, and all summands are regular outside $GP$.
Thus $Q$ is the only pole of $z$, and $e\,v_Q(z)=v_P(z)=-\ell$.

The $S$-orbit of $t$ has $e/|G_t|$ elements, each with pole divisor
$mP$. Hence $A_t\in E^S$ has pole order $\ell=em/|G_t|$ at $P$ and
is regular elsewhere. The trace therefore has pole order $m/|G_t|$
at $Q$, which also gives $|G_t|\mid m$.
\end{proof}

For the Hermitian coordinates $x,y$, put $e=|S|$. The products
$A_x,A_y$ have pole divisors $eaP_\infty,ebP_\infty$ in $H$,
and Lemma~\ref{hhh:orbit-trace} gives
\begin{equation}\label{hhh:orbit-poles}
 (U)_\infty=aQ_\infty,\qquad (V)_\infty=bQ_\infty.
\end{equation}
Since $a\mid q$ and $b\mid q+1$, we have $\gcd(a,b)=1$.
\begin{remark}\label{hhh:full-products}
If $S=1$, the functions in \eqref{hhh:orbit-generators} are the
ordinary traces of $x$ and $y$ from $H$ to $F$.
If $G=S$, then $U=A_x$ and $V=A_y$, with no further trace.
The full product $\prod_{s\in S}s(t)=A_t^{|G_t|}$ would give
pole orders $q,q+1$ in $F$ after taking traces.
\end{remark}

\subsection{The equation and its coefficients}
\begin{lemma}\label{hhh:coprime-plane}
Let $E/k$ be a function field with perfect full constant field $k$.
Suppose $u,v\in E$ have pole divisors $aQ_0,bQ_0$ at a rational
place $Q_0$, with $a,b>0$ coprime. Then the following hold.
\begin{enumerate}
\item We have $E=k(u,v)$, and there is a unique polynomial $\Phi$
of the form
\begin{equation}\label{hhh:rr-weighted-equation}
 \Phi(T_1,T_2)=T_2^a-cT_1^b+
 \sum_{\substack{i,j\geq0\\ ai+bj<ab}}c_{ij}T_1^iT_2^j,
 \qquad c\in k^*,\quad c_{ij}\in k,
\end{equation}
such that $\Phi(u,v)=0$.
\item The polynomial $\Phi$ is absolutely irreducible of total
degree $\max(a,b)$, and every polynomial relation between $u$ and $v$
is a multiple of $\Phi$.
\end{enumerate}
\end{lemma}
\begin{proof}
By the pole-degree formula \cite[Theorem~1.4.11]{Stichtenoth2009},
$[E:k(u)]=a$ and $[E:k(v)]=b$. Hence $[E:k(u,v)]$ divides
$\gcd(a,b)=1$, so $E=k(u,v)$. Since $v$ is regular away from the only
pole of $u$, it is integral over $k[u]$
\cite[Theorem~3.2.6(b)]{Stichtenoth2009}. As $k[u]$ is integrally
closed, \cite[Proposition~3.3.1]{Stichtenoth2009} gives the monic
minimal polynomial of $v$ over $k(u)$ in the form
\[
 \Phi(u,T_2)=T_2^a+\sum_{j=0}^{a-1}p_j(u)T_2^j,\qquad p_j\in k[T_1].
\]
The nonzero terms $p_j(u)v^j$ have pole orders $a\deg p_j+bj$,
pairwise distinct modulo $a$. Their sum is $-v^a$, so the largest
of these orders is $ab$, attained only for $j=0$. Thus $\deg p_0=b$, and each monomial $T_1^iT_2^j$ other than
$T_2^a$ and $T_1^b$ satisfies $ai+bj<ab$.
This gives \eqref{hhh:rr-weighted-equation}, uniquely since it is monic
of minimal degree in $T_2$.
Division by the monic polynomial $\Phi$ in $T_2$ shows that it generates
the ideal of relations. Since $ai+bj\le ab$ implies $i+j\le\max(a,b)$,
the terms $T_2^a$ and $-cT_1^b$ give $\deg\Phi=\max(a,b)$.
By \cite[Corollary~3.6.8]{Stichtenoth2009}, $\Phi(u,T_2)$ is irreducible
over $\bar k(u)$ since $k$ is the full constant field of $E=k(u,v)$.
Gauss's lemma now gives absolute irreducibility of $\Phi$.
\end{proof}

For $a,b>1$, plane curves with equations of the form
\eqref{hhh:rr-weighted-equation} are called $C_{ab}$ curves, up to
interchanging the coordinates \cite[Definition~1.1]{KawakitaMiura2002}.
Since $k$ is algebraically closed in $H$, it is also algebraically closed in $F$.
Applying the lemma to \eqref{hhh:orbit-poles} gives $F=k(U,V)$ and
a plane equation of degree $\max(a,b)\leq q+1$. Every irreducible
relation is a scalar multiple of $\Phi$, and monicity in $T_2$ fixes the scalar.

The weight of $U^iV^j$ is its pole order $ai+bj$ at $Q_\infty$.
By Lemma~\ref{hhh:coprime-plane}, the defining relation is supported
on the monomials indexed by $I$, where
\begin{equation}\label{hhh:rr-weighted-monomials}
 I=\{(i,j)\in\mathbb Z_{\ge0}^2:ai+bj\le ab\},\qquad
 \mathcal W=\{ai+bj:(i,j)\in I\}.
\end{equation}
To recover the relation, write it as
$\sum_{(i,j)\in I}d_{ij}U^iV^j=0$ with unknown coefficients
$d_{ij}\in k$, and expand at $P_\infty$. The pole order of $U^iV^j$ there is $e(ai+bj)$;
the following corollary shows that the coefficients of $\pi^{-ew}$,
$w\in\mathcal W$, already determine $(d_{ij})$ up to a scalar.
Let $L$ be the matrix with rows indexed by $\mathcal W$ and
columns indexed by $I$, with entries
\begin{equation}\label{hhh:laurent-matrix}
         L_{w,(i,j)}=[\pi^{-ew}]\bigl(U^iV^j\bigr),
\end{equation}
where $[\pi^n]f$ denotes the coefficient of $\pi^n$ in the Laurent
expansion of $f$ at $P_\infty$ with respect to $\pi$.

\begin{corollary}\label{hhh:laurent-elimination}
The matrix $L$ in \eqref{hhh:laurent-matrix} has a one-dimensional
kernel over $k$. If $(d_{ij})\in\ker L$ is normalized by $d_{0,a}=1$, then
$\Phi(T_1,T_2)=\sum_{(i,j)\in I}d_{ij}T_1^iT_2^j$ is the defining
polynomial in \eqref{hhh:rr-weighted-equation}.
The entries of $L$ are determined by the expansions of $U,V$ modulo
$\pi^{eab}k[[\pi]]$.
\end{corollary}
\begin{proof}
The pole order of $U^iV^j$ at $P_\infty$ is $e(ai+bj)$. Hence
\[
 L_{w,(i,j)}=0\quad(w>ai+bj),\qquad L_{ai+bj,(i,j)}\ne0.
\]
Two points of $I$ with the same weight differ by a multiple of $(b,-a)$.
Thus only $(b,0)$ and $(0,a)$ have the same weight. Reflection
$(i,j)\mapsto(b-i,a-j)$ pairs the lattice points of weight below
and above $ab$ in $[0,b]\times[0,a]$, so $|\mathcal W|=|I|-1=(a+1)(b+1)/2$.
Deleting column $(0,a)$ and ordering rows and columns by decreasing
weight gives a lower triangular square matrix with nonzero diagonal.
Thus $\dim_k\ker L=1$. Since $\Phi(U,V)=0$, the coefficient vector of
$\Phi$ lies in $\ker L$. Fixing $d_{0,a}=1$, we recover the other
coefficients by forward substitution in this triangular system.

To determine the required precision, consider a nonconstant monomial $U^iV^j$
of weight $m=ai+bj\le ab$.
Perturbing $U$ or $V$ by an element of $\pi^{eab}k[[\pi]]$ preserves
its pole order. Changing one factor $U$ in $U^iV^j$ in this way
introduces an error of valuation at least
\[
                         eab-e(m-a)\ge ea>0.
\]
For a factor $V$, the corresponding bound is $eab-e(m-b)\ge eb>0$.
Changing the factors one at a time shows that all coefficients
$[\pi^{-ew}](U^iV^j)$ with $w\ge0$ that enter $L$ are unchanged.
\end{proof}

Equation~\eqref{hhh:orbit-poles}, Lemma~\ref{hhh:coprime-plane} and
Corollary~\ref{hhh:laurent-elimination} prove Theorem~\ref{hhh:main}.

\begin{remark}\label{hhh:general-fields}
In the setting of Lemma~\ref{hhh:orbit-trace}, let $u,v\in E$ have pole
divisors $mP,nP$ with $m,n>0$ coprime. Applying the lemma to $u,v$
gives functions $U,V\in E^G$ with pole divisors $aQ,bQ$, where
$a=m/|G_u|$ and $b=n/|G_v|$. Since $a,b$ are coprime,
Lemma~\ref{hhh:coprime-plane} gives $E^G=k(U,V)$ and a relation of the
form \eqref{hhh:rr-weighted-equation}. The triangular system of
Corollary~\ref{hhh:laurent-elimination} recovers its coefficients whenever
the action of $G$ and local expansions at $P$ are computable.
If $g(E)\ge1$, Riemann--Roch allows the choice
$(m,n)=(2g(E),2g(E)+1)$ at every rational place. If $k=\F_{q^2}$
and $E/k$ is maximal, one may take $(m,n)=(q,q+1)$
\cite[Proposition~1.5(iv)]{FuhrmannGarciaTorres1997}.
\end{remark}

\subsection{Explicit computation}
For the required Laurent expansions, the equation
$\eta+\eta^q=\pi^{q+1}$ gives
\begin{equation}\label{hhh:laurent-coordinate}
 \eta=\sum_{r\geq0}(-1)^r\pi^{(q+1)q^r},\qquad
 B:=\pi^q x=\pi^{q+1}y=1+\eta^{q-1}.
\end{equation}
Indeed, adding the $q$-th power of the series to itself yields
$\pi^{q+1}$. Two solutions in $\pi k[[\pi]]$ differ by an element $d$
with $d(1+d^{q-1})=0$. The second factor is a unit, so the series is
the expansion of $Z/Y$ at $P_\infty$.
The identity for $B$ follows from $x=\pi/\eta$ and $y=1/\eta$.

\smallskip
\noindent\textbf{Input.} Matrix generators of $G$ acting on
columns $(X,Y,Z)^{\mathsf T}$ by left multiplication.
For $g\in G$, let $M_g$ be a matrix representing $g$.

\noindent\textbf{Output.} Rational expressions for $U,V$ in $x,y$ and their
defining polynomial $\Phi$ in the form \eqref{hhh:rr-weighted-equation}.

\begin{enumerate}
\item Enumerate $G$ and select the elements fixing $P_\infty$ to find $S$.
For $s\in S$, the third entry of $M_s^{-1}(x,y,1)^{\mathsf T}$ is a nonzero
scalar since $s$ preserves the tangent line $Z=0$ at $P_\infty$. Divide the
first two entries by the third to obtain $s(x),s(y)$ as in \eqref{hhh:point-stabilizer}.
For each coordinate, compare coefficients to identify the distinct affine-linear
forms, then multiply them to obtain $A_x,A_y$ in $k[x,y]$.
\item For each $P\in GP_\infty$, choose $g_P\in G$ with
$g_P(P_\infty)=P$, taking $g_{P_\infty}=1$. Put
$(\ell_0,\ell_1,\ell_2)^{\mathsf T}=M_{g_P}^{-1}(\pi,1,\eta)^{\mathsf T}$.
Since $g_P(f)=f\circ g_P^{-1}$, we have
\[
  g_P(A_t)=A_t(\ell_0/\ell_2,\ell_1/\ell_2) \qquad(t=x,y).
\]
Summing over $P$ and substituting $\pi=x/y$, $\eta=1/y$ gives $U,V$
as rational functions of $x,y$. For $P\ne P_\infty$, the point
$g_P^{-1}(P_\infty)$ is affine since $\mathcal H_q$ meets $Z=0$ only
at $P_\infty$. Hence $\ell_2$ has nonzero constant term and is a unit in $k[[\pi]]$.
\item Put $N=eab$ and $N_\ast=N+e\max(a,b)$. In \eqref{hhh:laurent-coordinate},
retain the terms of $\eta$ with $(q+1)q^r<N_\ast$ and compute
$B=1+\eta^{q-1}$ modulo $\pi^{N_\ast}$. For $P\ne P_\infty$, evaluate
the substitutions in (ii) modulo $\pi^N$. For $P=P_\infty$, use
\eqref{hhh:point-stabilizer} to write
\[
 \pi^q s(x)=\alpha B+\beta\pi^q,\qquad
 \pi^{q+1}s(y)=\alpha^{q+1}B+\alpha\beta^q\pi B+\gamma\pi^{q+1}.
\]
For each coordinate, multiply the distinct factors modulo $\pi^{N_\ast}$
to obtain $\pi^{ea}A_x$ and $\pi^{eb}A_y$, respectively. Since
$N_\ast-ea,N_\ast-eb\ge N$, dividing by $\pi^{ea}$ and $\pi^{eb}$ gives $A_x,A_y$
with errors in $\pi^N k[[\pi]]$. Sum over $P$ as in
\eqref{hhh:orbit-generators} to obtain $U,V$ modulo $\pi^N k[[\pi]]$.
\item Form $L$ by \eqref{hhh:laurent-matrix} and recover $\Phi$ by forward
substitution (Corollary~\ref{hhh:laurent-elimination}).
\end{enumerate}

\section{Genus bound and affine smoothness}\label{hhh:geometric}

To prove Theorem~\ref{hhh:defect}, we compute the Weierstrass semigroup
at $P_\infty\cap H^D$ for $D=G_xG_y$. We first recall the formula
of Campillo and Farr\'an \cite[Proposition~3.1]{CampilloFarran2000}
for the total affine delta invariant in terms of pole orders. For two
coprime pole orders it gives the genus bound and the affine smoothness
criterion for $C_{ab}$ curves \cite[Theorem~1.1(ii)]{KawakitaMiura2002}.

\begin{lemma}\label{hhh:pole-normalization}
Let $E/k$ be a function field with perfect full constant field $k$.
Suppose $u,v\in E$ have pole divisors $aQ_0,bQ_0$ at a rational place
$Q_0$, with $a,b>0$ coprime. Put $R=k[u,v]$, let $\widetilde R$ be its
integral closure in $E$, and set $\delta_{\mathrm{aff}}=\dim_k(\widetilde R/R)$.
Writing $\Gamma(Q_0)$ for the Weierstrass semigroup at $Q_0$, we have
\[
 \delta_{\mathrm{aff}}
 =\#\bigl(\Gamma(Q_0)\setminus\langle a,b\rangle\bigr)
 =\frac{(a-1)(b-1)}2-g(E).
\]
The affine model is smooth if and only if $\Gamma(Q_0)=\langle a,b\rangle$.
\end{lemma}
\begin{proof}
The ring $R$ lies in $\mathcal O_Q$ precisely for $Q\ne Q_0$, so
\cite[Theorem~3.2.6(b)]{Stichtenoth2009} gives
$\widetilde R=\bigcap_{Q\ne Q_0}\mathcal O_Q
=\bigcup_{n\ge0}\mathcal L(nQ_0)$.
The monic minimal polynomial in Lemma~\ref{hhh:coprime-plane} gives
$R=\bigoplus_{j=0}^{a-1}k[u]v^j$. The basis monomials $u^iv^j$
have distinct pole orders exhausting $\langle a,b\rangle$, so those of
pole order at most $n$ form a basis of $R\cap\mathcal L(nQ_0)$.
Since $Q_0$ is rational, $\dim_k\mathcal L(nQ_0)$ counts the
non-gaps at $Q_0$ up to $n$. Hence
\[
 \dim_k\frac{\mathcal L(nQ_0)}{R\cap\mathcal L(nQ_0)}
 =\#\bigl((\Gamma(Q_0)\setminus\langle a,b\rangle)\cap[0,n]\bigr).
\]
These quotients exhaust $\widetilde R/R$, which gives the first equality.
For large $n$, Riemann--Roch gives $\dim_k\mathcal L(nQ_0)=n+1-g(E)$,
so $\Gamma(Q_0)$ has $g(E)$ gaps. Subtracting this from Sylvester's
gap count $(a-1)(b-1)/2$ for $\langle a,b\rangle$ gives the second equality.
Finally, $\delta_{\mathrm{aff}}=0$ exactly when $R$ is integrally closed.
For the one-dimensional Noetherian domain $R$ this is equivalent to regularity,
and over the perfect field $k$ to smoothness.
\end{proof}

\begin{proposition}\label{hhh:semigroup}
Let $E/k$ be a function field with perfect full constant field $k$,
and let $G\le\operatorname{Aut}_k(E)$ be finite. Suppose $u,v\in E$
have pole divisors $mP,nP$ at a rational place $P$, where $1<m<n$ are
coprime and $\Gamma(P)=\langle m,n\rangle$. For $t=u,v$ let $G_t$ be
the stabilizer of $t$ in $G$, and put $a=m/|G_u|$, $b=n/|G_v|$.
Then $D=G_uG_v$ is a subgroup of $G$. Write $P_D=P\cap E^D$. For the
affine model of $F=E^G$ defined by \eqref{hhh:orbit-generators} using
$u,v$ and $S=\operatorname{Stab}_G(P)$, we have
\[
 \Gamma(P_D)=\langle a,b\rangle,\qquad
 \delta_{\mathrm{aff}}=g(E^D)-g(F).
\]
If $g(F)\ge2$, this model is smooth if and only if $G=D$.
\end{proposition}
\begin{proof}
Write $d_u=|G_u|$, $d_v=|G_v|$. Both stabilizers fix $P$, the only
pole of each of $u,v$.
Since $P$ is rational and $0,m$ are its only non-gaps up to $m$,
$\mathcal L(mP)=\langle1,u\rangle_k$. Each $\sigma\in G_v$ thus sends
$u$ to $\alpha u+\beta$ with $\alpha\in k^*$ and $\beta\in k$. Hence
$\sigma G_u\sigma^{-1}=G_{\sigma(u)}=G_{\alpha u+\beta}=G_u$, so $D$ is
a subgroup. As $E=k(u,v)$ by Lemma~\ref{hhh:coprime-plane},
$G_u\cap G_v=1$ and $|D|=d_ud_v$. The rational place $P$ is fixed by
$D$, hence totally ramified over $P_D$.

The stabilizers of $u,v$ in $D$ are still $G_u,G_v$, respectively.
By Lemma~\ref{hhh:orbit-trace} applied to $D$, the products over $D\cdot u$
and $D\cdot v$ have pole divisors $aP_D$ and $bP_D$ in $E^D$.
Thus $\langle a,b\rangle\subseteq\Gamma(P_D)$. Conversely, for
$r\in\Gamma(P_D)$, total ramification gives $d_ud_vr\in\Gamma(P)$.
Write $d_ud_vr=mi+nj$ with $i,j\ge0$. Since $d_u\mid m$, $d_v\mid n$
and $\gcd(m,n)=1$, reduction modulo $d_v$ and $d_u$ gives $d_v\mid i$
and $d_u\mid j$. Thus $r=a(i/d_v)+b(j/d_u)\in\langle a,b\rangle$.

Counting gaps at $P_D$ gives $g(E^D)=\tfrac12(a-1)(b-1)$.
Applying Lemma~\ref{hhh:orbit-trace} to $G$ gives $U,V$ with pole
divisors $aQ,bQ$ in $F$, where $Q=P\cap F$. Lemma~\ref{hhh:pole-normalization}
applied to $U,V$ in $F$ now gives $\delta_{\mathrm{aff}}=g(E^D)-g(F)$.
The model is smooth exactly when this difference vanishes. For $g(F)\ge2$
it is positive exactly when $D<G$, by the Hurwitz formula
\cite[Theorem~3.4.13]{Stichtenoth2009}
\[
 g(E^D)-g(F)=([G:D]-1)(g(F)-1)
                 +\tfrac12\deg\operatorname{Diff}(E^D/F).\qedhere
\]
\end{proof}

\begin{proof}[Proof of Theorem~\ref{hhh:defect}]
Apply Lemma~\ref{hhh:pole-normalization} to $U,V$ at $Q_\infty$ in $F$
to obtain the first equality of \eqref{hhh:defect-eq} and the criterion
for equality in the genus bound. Since $\Gamma(P_\infty)=\langle q,q+1\rangle$,
Proposition~\ref{hhh:semigroup} applied to $x,y$ at $P_\infty$ in $H$
gives $\delta_{\mathrm{aff}}=g(H^D)-g(F)$ and the last assertion.
\end{proof}

\smallskip
\noindent\emph{An explicit equation for $H^D$.} Put $d_y=|G_y|$. By
\eqref{hhh:point-stabilizer}, $G_y$ scales $x$ by the $d_y$th roots of unity
and $G_x$ translates $y$ by a subgroup $\Delta$ of
$\Lambda=\{c\in k:c^q+c=0\}$.
For a finite additive subgroup $J\subset k$, write
$B_J(T)=\prod_{c\in J}(T-c)$ for its separable additive polynomial
\cite[Lemma~3.1]{MaWangDecomposition}.
Put $x_D=x^{d_y}$ and $y_D=B_\Delta(y)$. The products over $D\cdot x$
and $D\cdot y$ are $(-1)^{d_y-1}x_D$ and $y_D$, respectively.
Since $D$ fixes $P_\infty$, Theorem~\ref{hhh:main} gives $H^D=k(x_D,y_D)$.

Put $M=B_{B_\Delta(\Lambda)}$. Since $\ker B_\Delta=\Delta$,
$\deg M=q/|\Delta|=a$. The monic polynomial $M\circ B_\Delta$ has
degree $q$ and vanishes on $\Lambda$, so it equals $T^q+T$. Thus
\[
                    x_D^b=M(y_D).
\]
This is the defining equation since $\deg M=[H^D:k(x_D)]$. The affine
model is smooth since the separable additive polynomial $M$ has nonzero
constant derivative.

\section{A family of \texorpdfstring{$S_3$}{S3}-quotients}\label{sec:s3}
Let $q=2^h$ with $h\ge1$, and let $G\simeq S_3$ act by coordinate
permutations on the Fermat model
$X_0^{q+1}+X_1^{q+1}+X_2^{q+1}=0$ of $\mathcal H_q$ over $k=\F_{q^2}$.
The only projective point fixed by $G$ is $(1:1:1)$, which lies off the curve.
Put $(\ell_0,\ell_1,\ell_2)=(X_1+X_2,X_2+X_0,X_0+X_1)$ and take
\[
 x=\frac{X_0+X_1+X_2}{\ell_2},\qquad y=\frac{\ell_1}{\ell_2}.
\]
The inverse formulas $(X_0:X_1:X_2)=(x+y+1:x+y:x+1)$ turn the Fermat
equation into $x^{q+1}=y^q+y$. Thus $x,y$ are Hermitian coordinates
with their only pole at $P_*=(1:1:0)$, the unique point of $\mathcal H_q$ on $\ell_2=0$.
Its stabilizer $S$ in $G$ is generated by the transposition interchanging
$X_0,X_1$, which fixes $x$ and sends $y$ to $y+1$. Hence $G_x=S$, $G_y=1$,
$A_x=x$ and $A_y=w:=y(y+1)$. Since $y^2+y=w$ and the nonidentity
element of $S$ moves $y$, we have $[H:k(x,w)]=2$ and $H^S=k(x,w)$.

The numerators of $x$ and $w=\ell_0\ell_1\ell_2/\ell_2^3$ are symmetric
in the Fermat coordinates. The summands in the traces of $A_x$ and $A_y$
in \eqref{hhh:orbit-generators} are therefore $x/d$ and $w/d^3$, respectively,
with $d$ running through $\ell_0/\ell_2=y+1$, $\ell_1/\ell_2=y$ and $1$.
Put $J=1+1/y+1/(y+1)=(w+1)/w$. Since $y^3+(y+1)^3=w+1$, the traces are
\begin{equation}\label{hhh:s3-generators}
 U=xJ,\qquad V=w+\frac{w+1}{w^2},\qquad
 V+1=wJ^3=\frac{(w+1)^3}{w^2}.
\end{equation}
Applying Theorem~\ref{hhh:main} in these Hermitian coordinates gives
$F=H^G=k(U,V)$, and the pole orders at
$Q_\infty=P_*\cap F$ are $a=q/2$ and $b=q+1$.

\begin{proposition}\label{hhh:s3-result}
Let $q=2^h$ with $h\ge1$, and let $G\simeq S_3$ permute the Fermat
coordinates of $\mathcal H_q$. With $U,V$ as in \eqref{hhh:s3-generators},
an affine plane model of the fixed field $F=H^G$ is given by
\begin{equation}\label{hhh:s3-equation}
                         U^{q+1}=\Psi_h(V+1),
\end{equation}
where the polynomials $\Psi_j\in\F_2[T]$ are defined by
\begin{equation}\label{hhh:s3-recurrence}
 \Psi_1=T,\qquad \Psi_2=T^2,\qquad
 \Psi_j=T^{2^{j-2}}(\Psi_{j-1}+\Psi_{j-2})\quad(j\ge3).
\end{equation}
Each $\Psi_j$ is monic of degree $2^{j-1}$. The model is smooth for $h=1$.
For $h\ge2$ its only geometric affine
singularity is $(U,V)=(0,1)$. For all $h\ge1$,
\[
 (\delta_{\mathrm{aff}},g(F))=
 \begin{cases}
 \left(\dfrac{(q+1)(q-2)}6,\ \dfrac{(q-2)^2}{12}\right), & h\text{ odd},\\[5pt]
 \left(\dfrac{q(q-1)}6,\ \dfrac{q(q-4)}{12}\right), & h\text{ even}.
 \end{cases}
\]
\end{proposition}
\begin{proof}
Put $\rho(w)=V+1=wJ^3=(w+1)J^2$, $L_j(w)=\sum_{r=0}^{j-1}w^{2^r}$
and $C_j=J^{2^j+1}L_j(w)$. Then $C_1=\rho(w)$ and $C_2=\rho(w)^2$.
For $j\ge3$, put $n=2^{j-2}$. Since $L_j=L_{j-1}+w^{2n}$ and
$L_{j-1}+L_{j-2}=w^n$, we obtain
\[
 L_j=(w+1)^nL_{j-1}+w^nL_{j-2},\qquad
 C_j=\rho(w)^n(C_{j-1}+C_{j-2}).
\]
Thus $C_j=\Psi_j(\rho(w))$ by \eqref{hhh:s3-recurrence}. For $j=h$,
$L_h(w)=y^q+y=x^{q+1}$, giving \eqref{hhh:s3-equation}. By the recurrence
$\Psi_h$ is monic of degree $q/2=[F:k(U)]$, so this is the defining equation.

The partial derivatives of the defining polynomial in \eqref{hhh:s3-equation} are $U^q$ and
$\Psi_h'(V+1)$. Hence its geometric affine singularities are precisely the
points $(0,V_0)$ for which $V_0+1$ is a multiple root of $\Psi_h$.
Let $m_h=\operatorname{ord}_{T=0}\Psi_h$. At $w=1$, $\rho$ has a triple
zero and $J$ a simple zero. Taking orders in
$\Psi_h(\rho(w))=J^{q+1}L_h(w)$ gives
\[
 3m_h=q+1+\operatorname{ord}_{w=1}L_h(w)=
 \begin{cases}q+1,&h\text{ odd},\\q+2,&h\text{ even},\end{cases}
\]
since $L_h'=1$ and $L_h(1)=h\bmod 2$. If $T_0\ne0$ is a root of $\Psi_h$,
choose $w_0\in\mathbb P^1(\bar k)$ with $\rho(w_0)=T_0$. Then $w_0\ne0,1,\infty$,
so $J$ is a unit at $w_0$. As $L_h'=1$, $\Psi_h(\rho(w))$ has a simple zero
at $w_0$, and hence $T_0$ is simple. For $h=1$ the equation is $V=U^3+1$.
For $h\ge2$ we have $m_h\ge2$, and $\delta_{\mathrm{aff}}$ equals the local
delta invariant at the unique affine singularity $(0,1)$
\cite[Chapter~IV, Exercise~1.8]{Hartshorne1977}.

Put $T=V+1$. Since $q+1$ is odd, the unit $\Psi_h(T)/T^{m_h}$ has a
$(q+1)$st root $\varepsilon$ in $\bar k[[T]]$. With $U_1=U/\varepsilon$
the completed local equation is $U_1^{q+1}=T^{m_h}$. For odd $h$,
$3m_h=q+1$ is odd, so the $m_h$ branches $T=\omega U_1^3$,
$\omega^{m_h}=1$, are distinct and smooth. Each pair has intersection
multiplicity $3$, giving $\delta_{\mathrm{aff}}=3\binom{m_h}{2}$
by \cite[Proposition~2.13]{BeelenPellikaan2000}.

For even $h$, $q+1=3m_h-1$ is coprime to $m_h$. The single branch has
completed local ring $\bar k[[\tau^{m_h},\tau^{q+1}]]$ with normalization
$\bar k[[\tau]]$. The gaps of $\langle m_h,q+1\rangle$ therefore give
$\delta_{\mathrm{aff}}=q(m_h-1)/2$.

Finally, $D=G_xG_y=S$, so
Theorem~\ref{hhh:defect} gives $g(F)=g(H^S)-\delta_{\mathrm{aff}}$
with $g(H^S)=(a-1)(b-1)/2=q(q-2)/4$. Substituting $m_h=(q+1)/3$ for
odd $h$ and $m_h=(q+2)/3$ for even $h$ gives the stated values of
$\delta_{\mathrm{aff}}$ and $g(F)$.
\end{proof}

Setting $a=e=1$ in the notation of
\cite[arXiv version, Proposition~3.7(ii), equation~(11)]{DallaVoltaMontanucciZini2019}
gives the same genera. For $q=8$ the equation is
$U^9=V(V+1)^3$ and $(\delta_{\mathrm{aff}},g(F))=(9,3)$.

\section*{Declaration of AI use}
The authors devised the main construction.
ChatGPT (GPT 6 Astra) was used for language editing, literature searches, candidate
constructions and proofs, and symbolic computations.
Claude (Fable 5) was used for language polishing and suggested the
analysis of the affine singularities in Section~\ref{sec:s3}.
The authors take full responsibility for all content of the manuscript.

\end{document}